\documentclass{aa}  

\usepackage{natbib}
\bibpunct{(}{)}{;}{a}{}{,} 

\usepackage{graphicx}

\usepackage{amssymb,amsmath}
\usepackage{delarray}
\usepackage{mathrsfs}
\usepackage{xcolor}
\usepackage[varg]{txfonts}
\usepackage{hyperref}
\hypersetup{
  colorlinks   = true, 
  urlcolor     = blue, 
  linkcolor    = blue, 
  citecolor   = blue 
}

\begin{document} 


\title{A novel fourth-order WENO interpolation technique}
\subtitle{A possible new tool designed for radiative transfer}

\author{Gioele Janett\inst{1,}\inst{2}
       \and
       Oskar Steiner\inst{1,}\inst{3}
       \and
       Ernest Alsina Ballester\inst{1}
       \and
       Luca Belluzzi\inst{1,}\inst{3}
       \and
       Siddhartha Mishra\inst{2}}
\institute{Istituto Ricerche Solari Locarno (IRSOL), 6605 Locarno-Monti, Switzerland
          \and
          Seminar for Applied Mathematics (SAM) ETHZ, D-MATH, 8093 Zurich, Switzerland
          \and
          Kiepenheuer-Institut f\"ur Sonnenphysik (KIS),
          Sch\"oneckstrasse 6, 79104 Freiburg i.~Br., Germany
\\          \email{gioele.janett@irsol.ch
          }
          }
 
\abstract
{
Several numerical problems require the interpolation of discrete data that
present at the same time (i) complex smooth structures and (ii) various types
of discontinuities. The radiative transfer
in
solar and stellar atmospheres is a typical example of such a problem. This calls for high-order well-behaved techniques that are able to
interpolate both smooth and discontinuous data.}
{This article expands on different nonlinear interpolation techniques
capable of guaranteeing
high-order accuracy and handling discontinuities
in an accurate and non-oscillatory fashion. {The final aim is to propose
new techniques which could be suitable for applications in the context of
numerical radiative transfer.}
}
{ We have proposed and tested two different techniques.
Essentially non-oscillatory (ENO) techniques generate several candidate interpolations based on different substencils.
The smoothest candidate interpolation is determined from a measure for the local smoothness,
thereby enabling the essential non-oscillatory property.
Weighted ENO (WENO) techniques {use} a convex
combination of all candidate substencils to obtain high-order accuracy in smooth regions while keeping the 
essentially non-oscillatory property.
In particular, we have outlined and tested
a novel well-performing fourth-order WENO interpolation technique for both uniform and nonuniform grids.}
{ { Numerical tests prove that the fourth-order WENO interpolation guarantees
fourth-order accuracy in smooth regions of the interpolated functions.
In the presence of discontinuities, {the fourth-order WENO interpolation} enables the non-oscillatory property, avoiding oscillations.
Unlike B\'ezier and monotonic {high-order} Hermite
interpolations, it does not degenerate to a linear interpolation near smooth extrema of the interpolated function.}}
{The novel fourth-order WENO interpolation
guarantees high accuracy in smooth regions, while effectively handling discontinuities.
This interpolation technique might be particularly suitable for several problems,
including a number of radiative transfer applications such as multidimensional problems,
multigrid methods, and formal solutions.}


\maketitle

\section{Introduction}\label{sec:sec1}
%
{In astrophysics, it is common practice
to solve the radiative transfer equation in solar and stellar atmospheres by means of numerical methods.
It is also known that realistic atmospheric models can be highly inhomogeneous and dynamic and often present discontinuities or sharp gradients.
Interpolations are ubiquitous in such problems:
for instance, they are a fundamental ingredient of multidimensional radiative transfer,
of multigrid methods, and of many formal solvers.
{For such problems, which contain both sharp} discontinuities
and complex smooth features,
one seeks interpolation techniques able to handle singularities
in an accurate and non-oscillatory fashion and at the same time
to perform with uniform high-order
accuracy in smooth regions.

It is a well known fact that standard high-order interpolations
tend to misrepresent the nonsmooth regions of
a problem, introducing spurious oscillations near discontinuities \citep[e.g.,][]{richards1991,zhang1997,shu1998}.
However, monotonicity-preserving strategies such as B\'ezier and monotonic high-order Hermite
interpolations, usually sacrifice accuracy to obtain monotonicity \citep{arandiga2013}. 
In fact, near smooth extrema of the interpolated function, they degenerate to a linear interpolation,
dropping to second-order accuracy \citep{shu1998}.
%
By contrast,} essentially non-oscillatory (ENO) and weighted ENO (WENO)
techniques accomplish the feat of approximating a function with high-order
accuracy in smooth parts, while avoiding “Gibbs-like” oscillations near the
discontinuities \citep{fjordholm2016}.

Section~\ref{sec:eno_approx} presents the ENO approximation technique
and it exposes the third-order ENO interpolation.
Section~\ref{sec:weno_approx} presents the WENO approximation technique,
with a special focus on a novel fourth-order WENO interpolation.
Section~\ref{sec:numerical_tests} exposes numerical tests for different interpolations.
Section~\ref{sec:RT_applications} describes different radiative transfer problems
that require high-order robust interpolation techniques.
Finally, Section~\ref{sec:conclusions} provides remarks and conclusions.
\section{ENO approximation}\label{sec:eno_approx}
The term stencil indicates the arrangement of grid points and corresponding function values (or coefficients) that is considered
for a particular interpolation or reconstruction.
Fixed stencil interpolations of second or higher order of accuracy are necessarily oscillatory
near a discontinuity and such oscillations do not decay in magnitude when the computational grid is refined; 
an explicit example {will be given in Section~\ref{sec:numerical_tests}.}
Such issues are well recognized and
considerable efforts have been already exercised to resolve them.
In a seminal paper, \citet{vanleer1979} proposed a hybrid technique that switches from a linear polynomial (second-order accurate)
to a piecewise constant approximation (first-order accurate) near discontinuities.
Similarly, one can also make use of a quadratic polynomial interpolation and limit the polynomial degree near discontinuities.
Alternatively, 
monotonic Hermite interpolants suppress over- and undershoots by controlling the first-order derivatives \citep{fritsch1980}, whereas
B{\'e}zier interpolants make use of the so-called control points to shape
the curve and suppress spurious extrema, preserving monotonicity in the interpolation.
These limiting strategies are very robust, produce smooth solutions, and are widely used for problems containing shocks and other discontinuities.
However,  
{such monotonicity-preserving strategies sacrifice accuracy in order to obtain monotonicity \citep{arandiga2013}. Their} main weakness is that
they necessarily degenerate
to {a linear interpolation} near smooth extrema of the interpolated function \citep{shu1998}.
{A local extremum is
any point at which the value of a function is larger (a maximum) or smaller (a minimum) than all the adjacent function values.
In this context, a smooth extremum is a local extremum of the function which is at least of class $C^2$, that is, twice continuously differentiable.}

By contrast, ENO techniques
are {high-order approximations}
that handle discontinuities and maintain high-order accuracy near smooth extrema.
The key point of the ENO strategy is to adaptively choose the stencil
in the direction of smoothness.
In practice, ENO schemes generate several candidate interpolations, measure the local smoothness
and choose the smoothest candidate interpolation to work with, discarding the rest.
This strategy enables the essential non-oscillatory property.
There is broad and well-cited body of  literature attesting to the importance
and the success of this strategy:
from the pioneering works by \citet{harten1986,harten1987} and~\citet{shu1988a,shu1988b}
to the recent survey by \citet{zhang2016}.

Most of the literature about ENO (and WENO) techniques is dedicated to reconstruction schemes and not to interpolations.
However, the reconstruction of cell averages of a function is equivalent to the interpolation of the
point values of its primitive. Therefore, ENO (and WENO) algorithms can be formulated
as an interpolation technique \citep{fjordholm2011}.
%
%
\subsection{Divided differences}\label{sec:div_diff}
Given a set of function values $\{y_k\}\coloneqq \{y(x_k)\}$ at positions $\{x_k\}$,
the zeroth degree divided differences are defined as
\begin{equation*}
y[x_k] \coloneqq y_k\,.
\end{equation*}
while the first degree divided differences read
\begin{equation*}
y[x_k,x_{k+1}]= \frac{y_{k+1}-y_k}{x_{k+1}-x_k}\,.
\end{equation*}
In general, the $j$-th degree forward divided differences, for $j>1$,
are inductively defined by
\begin{equation*}
y[x_k,\dots,x_{k+j}] \coloneqq \frac{y[x_{k+1},\dots,x_{k+j}]-y[x_k,\dots,x_{k+j-1}]}{x_{k+j}-x_k}\,.
\end{equation*}
%
As long as the function
$y(x)$ is smooth inside $[x_k,x_{k+j}]$, the
standard approximation theory affirms that \citep[e.g.,][]{fjordholm2016}
\begin{equation*}
y[x_k,\dots,x_{k+j}] = \frac{y^{(j)}(\xi)}{j!}\,,\text{ for some } \xi \in [x_{k},x_{k+j}]\,.
\end{equation*}
By contrast, if the function $y(x)$ is discontinuous at some point inside
$[x_k,x_{k+j}]$, then\footnote{For a local uniform grid, i.e., $x_{i+1}-x_i = \Delta x \;\forall\; i\in\{k,\dots,k+j-1\}$.}
\begin{equation*}
y[x_k,\dots,x_{k+j}] = \mathcal O\left(\frac{1}{\Delta x^j}\right)\,,
\end{equation*}
which diverges for $\Delta x\rightarrow 0$.
Thus, divided differences are a good measure of the smoothness
of the function $y(x)$ inside the stencil.
This property plays an essential role in the adaptive choice of the stencil in ENO schemes.
\subsection{ENO interpolation}\label{sec:eno_inter}
Consider a function $y(x)$ and
a set of function values $\{y_i\}$ at positions $\{x_i\}$ with $i\in\{1,\dots,N\}$.
The general $p$-th order ENO interpolation $E(x)$ in the interval $[x_1,x_N]$ reads
\begin{equation}\label{eno_inter}
E(x)=\sum_{i=1}^{N-1} p_i(x){I_i(x)}\,,
\end{equation}
with
\begin{equation}\label{pw_constant_I}
I_i(x)=
\begin{cases}
1 & \text{ if } x\in[x_i,x_{i+1}]\,, \\ 
0 & \text{ if } x\notin[x_i,x_{i+1}] \,,
\end{cases}
\end{equation}
and $p_i(x)$ denotes the interpolating polynomial acting on the interval $[x_i,x_{i+1}]$.
This polynomial must satisfy the accuracy requirement
\begin{equation*}
p_i(x) = y(x) + \mathcal O (\Delta x^p)\,,\text{ for } x\in[x_i,x_{i+1}].
\end{equation*}
The choice of ENO stencil for constructing $p_i(x)$ is
based on the divided differences presented in Section~\ref{sec:div_diff}.
The stencil is adaptively chosen {so} that,
starting from $x_i$, it extends in the direction where $y(x)$ is smoothest,
{that is}, where the divided differences are smallest.

A general $p$-point stencil containing $x_i$ is given by
\begin{equation}\label{stencil}
S_p^q = \{x_{i-p+q},\dots,x_{i+q-1}\}\,,
\end{equation}
with $q\in\{1,\dots,p\}$. For notational simplicity,
the explicit dependence of the stencil $S_p^q$ on index $i$ is not indicated. 
This stencil can be expanded in two ways: by adding either the left neighbor $x_{i-p+q-1}$
or the right neighbor $x_{i+q}$.
One decides upon which point to add by comparing
the two relevant divided differences and picking the one with a smaller absolute value.
Thus, if
\begin{equation*}
|y[x_{i-p+q-1},\dots,x_{i+q-1}]|\le |y[x_{i-p+q},\dots,x_{i+q}]|\,,
\end{equation*}
one adds $x_{i-p+q-1}$ to the stencil, otherwise, $x_{i+q}$.

In order to build $p_i(x)$, one starts with the one-point stencil
\begin{equation*}
S_1^1 = \{x_i\}\,, 
\end{equation*}
and adaptively adds points with the procedure described above, ending up with a $p$-point stencil
\begin{equation*}
S_p^{q} = \{x_{i-p+q},\dots,x_{i+q-1}\}\,,\text{ for some }q\in\{1,\dots,p\}\,.
\end{equation*}
We note that each substencil $S_p^{q}$ is a substencil of the large stencil
\begin{equation*}
S_{2p-1}^p = \{x_{i-p+1},\dots,x_{i+p-1}\}\,, 
\end{equation*}
and always contains the considered $x_i$ (but not necessarily $x_{i+1}$).

Numerical analysis asserts that
there is a unique polynomial $q_q(x)$ of degree at most $p-1$,
{that} goes through all
the data points of the substencil $S_p^{q}$.
Its Lagrangian expression reads
\begin{equation}\label{lagrange_interp}
q_q(x)=\sum_{x_k\in S_p^{q}} y(x_k)\ell_k(x)\,,
\end{equation}
%
where the Lagrange basis polynomials $\ell_k$ are defined as\footnote{Lagrange basis polynomials
satisfy the relation $\ell_k(x_m)=\delta_{km}$, being $\delta_{km}$ the Kronecker delta given by
\begin{equation*}
\delta_{km} = 
\begin{cases}
1 & \text{ if } k = m \,, \\ 
0 & \text{ if } k \ne m\,.
\end{cases}
\end{equation*}}
\begin{equation*}
\ell_k(x)=\prod_{\substack{x_j\in S_p^{q}\\j\neq k}}\frac{x-x_j}{x_k-x_j}\,.
\end{equation*}
Hence, once the optimal candidate substencil $S_p^{q}$ is determined,
the $p$-order accurate polynomial $p_i(x)$ in Equation~\eqref{eno_inter} is simply chosen as
\begin{equation*}
p_i(x)=q_q(x)\,,
\end{equation*}
with $q_q(x)$ given by Equation~\eqref{lagrange_interp}.

Even though there are very {few} theoretical results about the stability of ENO
schemes \citep[e.g.,][]{fjordholm2012,fjordholm2011}, these schemes are very robust and stable in practice.
In fact, the ENO interpolant $E(x)$ is essentially non-oscillatory in the sense of recovering $y(x)$
to order $\mathcal O(\Delta x^p)$ near (but not at) the location of a discontinuity.
We note that ENO schemes
allow for minimal over- or undershoots near discontinuities. 
As an explicit example of this technique, the third-order accurate ENO interpolation is presented in the following.
\subsection{Third-order ENO interpolation}\label{sec:eno_inter3}
In order to build a third-order accurate $p_i(x)$, one needs a three-point stencil.
Thus, starting from the one-point stencil
$$S_1^1=\{x_i\}\,,$$
one adds either the left neighbor $x_{i-1}$
or the right neighbor $x_{i+1}$.
One takes the decision by comparing the absolute values of the two relevant divided differences
$$y [x_{i-1},x_i]=\frac{y_i-y_{i-1}}{x_{i}-x_{i-1}}\,,$$
$$y [x_i,x_{i+1}]=\frac{y_{i+1}-y_i}{x_{i+1}-x_i}\,.$$
The smaller one implies that the function is smoother in that stencil.
Therefore, if
$$|y [x_{i-1},x_i]| < |y [x_i,x_{i+1}]|\,,$$
the two-point stencil is taken as
\begin{equation*}
S^{1}_2 = \{x_{i-1},x_i\}\,, 
\end{equation*}
otherwise,
\begin{equation*}
S^{2}_2 = \{x_i,x_{i+1}\}\,. 
\end{equation*}
Suppose that the two-point stencil
$S^{2}_2=\{x_i,x_{i+1}\}$ has been chosen.
In this case, one can add either the left neighbor $x_{i-1}$
or the right neighbor $x_{i+2}$.
One takes the decision by comparing the absolute values of the two relevant divided differences
%
$$y [x_{i-1},x_i,x_{i+1}]=\frac{y[x_i,x_{i+1}]-y[x_{i-1},x_i]}{x_{i+1}-x_{i-1}}\,,$$
$$y [x_i,x_{i+1},x_{i+2}]=\frac{y[x_{i+1},x_{i+2}]-y[x_i,x_{i+1}]}{x_{i+2}-x_i}\,.$$
Once again, if
$$|y [x_{i-1},x_i,x_{i+1}]| < |y [x_i,x_{i+1},x_{i+2}]|\,,$$
the three-point stencil is taken as
\begin{equation*}
S^{2}_3 = \{x_{i-1},x_i,x_{i+1}\}\,, 
\end{equation*}
otherwise,
\begin{equation*}
S^{3}_3 = \{x_i,x_{i+1},x_{i+2}\}\,. 
\end{equation*}
Then, one proceeds with the Lagrange interpolation given by Equation~\eqref{lagrange_interp}
on the selected substencil, yielding the third-order accurate quadratic polynomial $p_i(x)$.

We note that the choice of the two-point stencil~$S^{1}_2 = \{x_{i-1},x_i\}$ at the previous step
may lead to the three-point stencil
\begin{equation*}
S^{1}_3 = \{x_{i-2},x_{i-1},x_i\}\,.
\end{equation*}
In order to build $p_i(x)$, the third-order ENO interpolation thus considers the large stencil
\begin{equation*}
S_5^3 = \{x_{i-2},x_{i-1},x_i,x_{i+1},x_{i+2}\}\,.
\end{equation*}
All the substencils $S^{1}_3$, $S^{2}_3$, and $S^{3}_3$ usually work well for globally smooth problems, while
the ENO technique adaptively avoids including possible discontinuities in the selected substencil.
%
%
\section{WENO approximation}\label{sec:weno_approx}
WENO approximations are inspired
by the ENO strategy, {offering} with additional advantages.
Instead of using the smoothest candidate substencil to form the interpolation, WENO schemes use a convex
combination of all candidate substencils to obtain high-order accuracy in smooth regions while keeping the 
essentially non-oscillatory property.
Additionally, the WENO strategy yields smoother interpolations and avoids the use
of conditional ``if$\dots$else'' statements, which
burden the algorithm.

\citet{liu1994} proposed this strategy, constructing the first WENO schemes.
\citet{jiang1996} constructed arbitrary-order accurate WENO
approximations for finite difference schemes and most applications use their fifth-order WENO version.
General information is found in the recent surveys by \citet{shu2009} and~\citet{zhang2016}.
Many papers about WENO approximations
consider their implementation to uniform grids.
However, the WENO strategy extends naturally to nonuniform grids, although it becomes quite complicated,
depending on the complexity
of the grid structure and on the interpolation degree \citep{crnjaric2007}.
\subsection{WENO interpolation}
Consider a function $y(x)$ and
a set of function values $\{y_i\}$ at positions $\{x_i\}$ with $i\in\{1,\dots,N\}$.
The general $p$-th order WENO interpolation $W(x)$ in the interval $[x_1,x_N]$ reads
\begin{equation}\label{weno_inter}
W(x)=\sum_{i=1}^{N-1} p_i(x){I_i(x)}\,,
\end{equation}
where $I_i(x)$ is defined in Equation~\eqref{pw_constant_I}.
Consider the $p$-point stencil $S_p^q$ given by Equation~\eqref{stencil}
to be the large stencil.
%
The unique $p$-th-order accurate Lagrange polynomial $p(x)$ interpolating $S_p^q$
can be recovered  as a linear combination of
$\ell$-th-order accurate Lagrange polynomials $q_m(x)$ defined on a set of $\ell$-point substencils $S_\ell^m$ with $\ell\in\{1,\dots,p\}$
and $m\in\{-p+q+\ell,\dots,q\}$, i.e.,\footnote{The cases with $m<1$ or $m>\ell$ are allowed. This means that some substencils $S_\ell^m$ may not contain $x_i$.}
\begin{equation}\label{linear_combination}
p(x)=\sum_{m=-p+q+\ell}^{q} \gamma_m(x)q_m(x)\,,
\end{equation}
where $\gamma_m(x)$ are the so-called linear (or optimal) weights,
which are always positive and univocally defined. Because of consistency one has
\begin{equation*}
\sum_{m=-p+q+\ell}^{q} \gamma_m(x)=1\,.
\end{equation*}
The key idea of the WENO strategy is to
build the polynomial $p_i(x)$ in Equation~\eqref{weno_inter}
as a weighted combination based on Equation~\eqref{linear_combination}, namely,
%
\begin{equation*}
p_i(x)=\sum_{m=-p+q+\ell}^{q} \omega_m(x)q_m(x)\,,
\end{equation*}
where the weights coefficients $\omega_m(x)$ must satisfy the requirements:
\begin{itemize}
 \item $\omega_m(x)\ge0\;  \forall\; m$ for stability;
 \item $\sum_{m=-p+q+\ell}^{q} \omega_m(x)=1$ for consistency;
 \item $\omega_m(x)\approx \gamma_m(x) \;  \forall\; m$ if $y(x)$ is smooth in the large stencil $S_p^q$;
 \item $\omega_m(x)\approx 0$ if $y(x)$ has a discontinuity in the stencil $S_\ell^m$.
\end{itemize}
These (nonlinear) weights are usually defined as
\begin{equation*}
\omega_m(x) = \frac{\alpha_m(x)}{\sum_{k=-p+q+\ell}^q \alpha_k(x)}\,,
\end{equation*}
where the unnormalized weights are given by
\begin{equation*}
\alpha_m(x) = \frac{\gamma_m(x)}{(\epsilon+\beta_m)^a}\,,
\end{equation*}
and $\epsilon$ is a small positive constant used to avoid vanishing denominators;
typically $\epsilon=10^{-6}$ \citep{zhang2016}.
{A} larger power $a$ makes the weight assigned to a nonsmooth substencil approach to zero faster,
resulting in more dissipative WENO interpolations \citep{liu2018}.
{The exponent $a$ is a free parameter. In this work, one has either $a=1$ or $a=3/2$.}

The nonlinear weights $\omega_m(x)$ rely on the smoothness indicators $\beta_m$,
which measure the relative smoothness of the function $y(x)$ inside the stencils $S_\ell^m$.
A larger $\beta_m$ indicates a lack of smoothness of $y(x)$ in the stencil $S_\ell^m$
and produces a smaller nonlinear weight $\omega_m(x)$.
The choice of the smoothness indicators plays an essential role
in the performance of WENO approximations.
Within the literature, the version of \citet{jiang1996}
\begin{equation}\label{smoothness_indicator}
\beta_m=\sum_{j=1}^{\ell-1} \Delta x^{2j-1} \int_{i}^{i+1}\left(\frac{{\rm d}^j}{{\rm d}x^j}q_m(x)\right)^2{\rm d}x\,,
\end{equation}
is well established.
Unfortunately, indicators~\eqref{smoothness_indicator} consider uniform grids alone
and are particularly dissipative (i.e., they smear sharp gradients) for low-order WENO interpolations.
In the following, a third-order and a novel fourth-order accurate WENO interpolations are presented.
\subsection{Third-order WENO interpolation}\label{sec:weno_inter3}
Compared with higher-order versions, the third-order WENO interpolation is more
robust for the treatment of discontinuous problems, it uses fewer grid points\footnote{Hence, it reduces the difficulty of the boundary treatment and
is easily generalized to unstructured grids.},
and it provides a suitable compromise between computational cost and 
accuracy \citep{liu2018}.

The unique quadratic Lagrange polynomial $p(x)$ interpolating the three-point large stencil
$$S_3^2 = \{x_{i-1}, x_i , x_{i+1}\}\,,$$
can be written as
$$p(x) = \gamma_1(x) q_1(x)+ \gamma_2(x) q_2(x)\,,$$
in other words, as a linear combination of the linear Lagrange interpolations
\begin{align*}
q_1(x) &= \frac{-y_{i-1}(x-x_i)+y_i(x-x_{i-1})}{x_i-x_{i-1}}\,,\\
q_2(x) &= \frac{-y_{i}(x-x_{i+1})+y_{i+1}(x-x_{i})}{x_{i+1}-x_{i}}\,,
\end{align*}
based, respectively, on the two two-point substencils
$$S_2^{1} = \{x_{i-1}, x_i\}\,,\; \text{ and }\; S_2^{2} = \{x_i , x_{i+1}\}\,.$$
The linear weights are given by
\begin{equation*}
\gamma_1(x) = -\frac{x-x_{i+1}}{x_{i+1}-x_{i-1}}\,,\;\text{ and }\;
\gamma_2(x) = \frac{x-x_{i-1}}{x_{i+1}-x_{i-1}}\,,
\end{equation*}
and satisfy $\gamma_1(x) + \gamma_2(x) = 1$.
The linear weights depend just on the grid
geometry and not on the function values.
The third-order WENO interpolation in the interval $[x_i,x_{i+1}]$ reads
$$p_i(x) = \omega_1(x) q_1(x)+ \omega_2(x) q_2(x)\,,$$
where the nonlinear weights are defined as
\begin{equation*}
\omega_1(x) = \frac{\alpha_1(x)}{\alpha_1(x)+\alpha_2(x)}\,,\;\text{ and }\;
\omega_2(x) = \frac{\alpha_2(x)}{\alpha_1(x)+\alpha_2(x)}\,,
\end{equation*}
with
\begin{equation*}
\alpha_1(x) = \frac{\gamma_1(x)}{(\epsilon+\beta_1)^\frac{3}{2}}\,,\;\text{ and }\;
\alpha_2(x) = \frac{\gamma_2(x)}{(\epsilon+\beta_2)^\frac{3}{2}}\,,
\end{equation*}
and $\epsilon=10^{-6}$.
\subsubsection{Smoothness indicators for uniform grids}
In uniform grids,
the conventional choice for the smoothness indicators~\eqref{smoothness_indicator} yields
\begin{equation*}
\beta_1 = (y_i-y_{i-1})^2\,,\;\text{ and }\; 
\beta_2 = (y_{i+1}-y_{i})^2\,. 
\end{equation*}
Unfortunately, this version renders the third-order WENO interpolation too dissipative.

Alternatively, the mapping function of the WENO-M scheme \citep{henrik2005} and
the global stencil indicator of the WENO-Z scheme \citep{borges2008} effectively improve the accuracy of high-order WENO schemes,
but they are not satisfactory for the third-order one.
More recently, \citet{wu2015} proposed the
global smoothness indicators for the less dissipative third-order WENO-N3 and WENO-NP3 schemes.
Similarly, \citet{gande2017} proposed the
global smoothness indicator of the WENO-F3 scheme.
However, such methods are usually not satisfactory in terms of accuracy or
generate apparent oscillatory interpolations.

For these reasons, \citet{liu2018} proposed smoothness indicators
{that make use} of all the three points of the $S_3^2$ stencil,
{namely,}
\begin{equation}\label{beta_liu}
\begin{aligned} 
\beta_1 &= \frac{1}{4}(|y_{i+1}-y_{i-1}|-|4y_i-3y_{i-1}-y_{i+1}|)^2\,,\\ 
\beta_2 &= \frac{1}{4}(|y_{i+1}-y_{i-1}|-|4y_i-y_{i-1}-3y_{i+1}|)^2\,. 
\end{aligned}
\end{equation}
Numerical tests demonstrate that these new
indicators provide less dissipation and better resolution than the standard ones.
\subsubsection{Smoothness indicators for nonuniform grids}
\citet{liu2018} defined the smoothness indicators~\eqref{beta_liu} for uniform grids only.
The generalization to the nonuniform case is presented in the following. For notational convenience
\begin{equation*}
d_i=\frac{y_{i+1}-y_i}{h_i}\,,\; \text{ with }\;h_i=x_{i+1}-x_i\,.
\end{equation*}
The first-order three-point numerical derivatives of $y(x)$ at nodes $x_{i-1}$, $x_i$ and $x_{i+1}$ for nonuniform grids are approximated by \citep[e.g.,][]{singh2009}
%
\begin{equation}\label{derivatives_nonuniform}
\begin{aligned}
y'_{i-1} &= \frac{2h_{i-1}+h_i}{h_{i-1}+h_i}\cdot d_{i-1}
-\frac{h_{i-1}}{h_{i-1}+h_i}\cdot d_i\,,\\ 
y'_i &= \frac{h_i}{h_{i-1}+h_i}\cdot d_{i-1}
+\frac{h_{i-1}}{h_{i-1}+h_i}\cdot d_i\,,\\ 
y'_{i+1} &= \frac{-h_i}{h_{i-1}+h_i}\cdot d_{i-1}
+\frac{h_{i-1}+2h_i}{h_{i-1}+h_i}\cdot d_i\,. 
\end{aligned}
\end{equation}
Following the underlying idea of \citet{liu2018}, the new indicators are constructed as follows
\begin{equation}\label{beta_liu_nonuniform}
\begin{aligned} 
\beta_1 &= h_i^2\left(|y'_{i}|-|y'_{i-1}|\right)^2\,,\\
\beta_2 &= h_{i-1}^2\left(|y'_{i+1}|-|y'_{i}|\right)^2\,.
\end{aligned}
\end{equation}
If $h_{i-1}=h_i$, the indicators~\eqref{beta_liu_nonuniform} reduce to the uniform case indicators~\eqref{beta_liu}.

In monotonic smooth regions,
the numerical derivatives~\eqref{derivatives_nonuniform} have the same sign, meaning that
$\operatorname{sign} y'_{i-1} = \operatorname{sign} y'_i = \operatorname{sign} y'_{i+1}$, and
the indicators~\eqref{beta_liu_nonuniform} reduce to
\begin{equation*}
\beta_1 = \beta_2  =
\left(2\frac{h_{i-1}h_i}{h_{i-1}+h_i}\right)^2\left( d_i- d_{i-1}\right)^2\,.
\end{equation*}
and, consequently,
\begin{equation*}
\omega_1(x)=\gamma_1(x)\,,\;\text{ and }\;\omega_2(x)=\gamma_2(x)\,.
\end{equation*}
In monotonic smooth regions, the nonlinear weights are exactly equal to the optimal weights,
reducing numerical dissipation for the WENO interpolation.

Assume now that $y(x)$ has {a} discontinuity in the interval $[x_{i-1},x_i]$
and is smooth in the interval $[x_i,x_{i+1}]$. In this case one has
\begin{equation*}
\left| d_{i-1}\right|\gg \left| d_i\right|\,.
\end{equation*}
Consequently, $\operatorname{sign} y'_{i-1} = \operatorname{sign} y'_i \ne \operatorname{sign} y'_{i+1}$ and
\begin{equation*}
\begin{aligned} 
\beta_1 &= \left(2\frac{h_{i-1}h_i}{h_{i-1}+h_i}\right)^2\left( d_i- d_{i-1}\right)^2=
\mathcal{O}\left( d_{i-1}^2\right)\,,\\
\beta_2 &= (2h_{i-1}d_i)^2=\mathcal{O}\left( d_i^2\right)\,,
\end{aligned}
\end{equation*}
resulting in $\beta_1\gg \beta_2$. 
Analogously, if $y(x)$ is smooth in $[x_{i-1},x_i]$
and has {a} discontinuity inside $[x_i,x_{i+1}]$, one obtains
\begin{equation*}
\begin{aligned} 
\beta_1 &=\mathcal{O}\left( d_{i-1}^2\right)\,,\\
\beta_2 &=\mathcal{O}\left( d_i^2\right)\,,
\end{aligned}
\end{equation*}
resulting in $\beta_1\ll \beta_2$.
This guarantees that the smoothness indicators~\eqref{beta_liu_nonuniform} effectively detect discontinuities in nonuniform grids and
{ensure} the non-oscillatory property in the interpolation.
The indicators~\eqref{beta_liu} and \eqref{beta_liu_nonuniform} do not
distinguish between discontinuities and smooth extrema.
For this reason, they do not guarantee accuracy
near such critical points.
Some numerical evidence for the third-order WENO interpolation with smoothness indicators~\eqref{beta_liu_nonuniform}
is shown in Figures~\ref{fig:exp}-\ref{fig:exp2_not}.
%
%
\subsection{Fourth-order WENO interpolation}\label{sec:weno_inter4}
{The WENO interpolation and the relative smoothness indicators
outlined in the following are completely original.
This novel fourth-order WENO interpolation} uses a four-point stencil,
is symmetric around the considered interval $[x_i , x_{i+1}]$,
and provides a suitable compromise of the computation cost and the accuracy.

The unique cubic Lagrange polynomial $p(x)$ interpolating the four-point large stencil
$$S_4^3 = \{x_{i-1}, x_i , x_{i+1}, x_{i+2}\}\,,$$
can be written as
$$p(x) = \gamma_2(x) q_2(x)+ \gamma_3(x) q_3(x)\,,$$
that is, as a linear combination of the quadratic Lagrange interpolations
{\small
\begin{align*}
q_2(x) &= y_{i-1}\frac{(x\!-\!x_i)(x\!-\!x_{i+1})}{h_{i-1}(h_{i-1}+h_i)}
\!-\!y_i\frac{(x\!-\!x_{i-1})(x\!-\!x_{i+1})}{h_{i-1}h_i}
\!+\!y_{i+1}\frac{(x\!-\!x_{i-1})(x\!-\!x_i)}{(h_{i-1}+h_i)h_i},\\
q_3(x) &= y_i\frac{(x\!-\!x_{i+1})(x\!-\!x_{i+2})}{h_{i}(h_{i}+h_{i+1})}
\!-\!y_{i+1}\frac{(x\!-\!x_i)(x\!-\!x_{i+2})}{h_ih_{i+1}}
\!+\!y_{i+2}\frac{(x\!-\!x_i)(x\!-\!x_{i+1})}{(h_i+h_{i+1})h_{i+1}},
\end{align*}}\noindent
based, respectively, on the two three-point substencils
$$S_3^{2} = \{x_{i-1}, x_i, x_{i+1}\}\,,\; \text{ and }\; S_3^{3} = \{x_i , x_{i+1}, x_{i+2}\}\,.$$
The linear weights are given by
\begin{equation*}
\gamma_2(x) = -\frac{x-x_{i+2}}{x_{i+2}-x_{i-1}}\,,\;\text{ and }\;
\gamma_3(x) = \frac{x-x_{i-1}}{x_{i+2}-x_{i-1}}\,,
\end{equation*}
and satisfy $\gamma_2(x) + \gamma_3(x) = 1$.
The linear weights depend just on the grid
geometry and not on the {values of the function}.
The fourth-order WENO interpolation in the interval $[x_i,x_{i+1}]$ reads
$$p_i(x) = \omega_2(x) q_2(x)+ \omega_3(x) q_3(x)\,,$$
where the nonlinear weights are defined as
\begin{equation*}
\omega_2(x) = \frac{\alpha_2(x)}{\alpha_2(x)+\alpha_3(x)}\,,\;\text{ and }\;
\omega_3(x) = \frac{\alpha_3(x)}{\alpha_2(x)+\alpha_3(x)}\,,
\end{equation*}
with
\begin{equation*}
\alpha_2(x) = \frac{\gamma_2(x)}{\epsilon+\beta_2}\,,\;\text{ and }\;
\alpha_3(x) = \frac{\gamma_3(x)}{\epsilon+\beta_3}\,,
\end{equation*}
and $\epsilon=10^{-6}$.
\subsubsection{Smoothness indicators for uniform grids}\label{subsec:uniform4}
The first-order 4-point numerical derivatives of $y(x)$ at nodes $x_{i-1}$, $x_i$, $x_{i+1}$
and $x_{i+2}$ for uniform grids read \citep[e.g.,][]{singh2009}
\begin{equation}\label{derivatives_uniform4}
      \begin{bmatrix}
      y'_{i-1}  \\
      y'_i      \\
      y'_{i+1}  \\
      y'_{i+2} 
      \end{bmatrix}=\frac{1}{6h}
      \begin{bmatrix}
      -11       & 18    & -9    & 2     \\
      -2        & -3    & 6     & -1    \\
      1         & -6    & 3     & 2     \\
      -2        & 9     & -18   & 11 
      \end{bmatrix}
      \begin{bmatrix}
      y_{i-1}  \\
      y_i       \\
      y_{i+1}   \\
      y_{i+2} 
      \end{bmatrix}\,,    
\end{equation}
or, alternatively,
\begin{equation*}
      \begin{bmatrix}
      y'_{i-1}  \\
      y'_i      \\
      y'_{i+1}  \\
      y'_{i+2} 
      \end{bmatrix}=\frac{1}{6}
      \begin{bmatrix}
      11        & -7    & 2     \\
      2         & 5     & -1    \\
      -1        & 5     & 2     \\
      2         & -7    & 11     
      \end{bmatrix}
      \begin{bmatrix}
      d_{i-1}  \\
      d_i       \\
      d_{i+1}    
      \end{bmatrix}\,.    
\end{equation*}
{Using these,} the new indicators are constructed as follows
\begin{equation}\label{beta_uniform4}
\begin{aligned} 
\beta_2 &= 4\left(|y'_{i+1}-y'_{i}|-|y'_{i}-y'_{i-1}|\right)^2\,,\\
\beta_3 &= 4\left(|y'_{i+2}-y'_{i+1}|-|y'_{i+1}-y'_{i}|\right)^2\,,
\end{aligned}
\end{equation}
with
\begin{equation}\label{derivatives_diff_uniform4}
      \begin{bmatrix}
      y'_i-y'_{i-1}     \\
      y'_{i+1}-y'_i     \\
      y'_{i+2}-y'_{i+1}
      \end{bmatrix}=\frac{1}{6}
      \begin{bmatrix}
      -9        & 12    & -3    \\
      -3        & 0     & 3     \\
      3         & -12   & 9      
      \end{bmatrix}
      \begin{bmatrix}
      d_{i-1}  \\
      d_i       \\
      d_{i+1}    
      \end{bmatrix}\,.    
\end{equation}

In smooth regions {where the second derivative has a constant sign}, the differences of
the numerical derivatives~\eqref{derivatives_diff_uniform4} have the same sign, that is,
$$\operatorname{sign} \left( y'_i -y'_{i-1}\right) = \operatorname{sign} \left( y'_{i+1} -y'_i\right) = \operatorname{sign} \left( y'_{i+2} -y'_{i+1}\right)\,.$$
Consequently,
the indicators~\eqref{beta_uniform4} reduce to
\begin{equation*}
\beta_2 = \beta_3 = 4\left(-y_{i-1}+3y_i-3y_{i+1}+y_{i+2}\right)^2\,,
\end{equation*}
and one gets
\begin{equation*}
\omega_2(x)=\gamma_2(x)\,,\;\text{ and }\;\omega_3(x)=\gamma_3(x)\,.
\end{equation*}
Hence, in smooth regions {where the second derivative has a constant sign}, the nonlinear weights correspond to the linear weights,
reducing the numerical dissipation of the fourth-order WENO interpolation.
We note that the indicators~\eqref{beta_uniform4} allow for over- or undershoots around smooth extrema,
yielding a higher accuracy
near such critical points with respect {to monotonic interpolations.
However, this feature may have the adverse effect of producing negative interpolated values in strictly
positive functions.
Moreover,} the indicators~\eqref{beta_uniform4} do not guarantee accuracy
near inflection points.

For the discontinuous case, assume that $y(x)$ has {a} discontinuity in the interval $[x_{i-1},x_i]$
and is smooth in the interval $[x_i,x_{i+2}]$. In this case one has
\begin{equation*}
\left| d_{i-1}\right|\gg \left| d_i\right|,\left| d_{i+1}\right|\,,
\end{equation*}
and from Equation~\eqref{derivatives_diff_uniform4}, one has
$$\operatorname{sign} \left( y'_i -y'_{i-1}\right) = \operatorname{sign} \left( y'_{i+1} -y'_i\right) \ne \operatorname{sign} \left( y'_{i+2} -y'_{i+1}\right)\,.$$
Consequently,
\begin{equation*}
\begin{aligned} 
\beta_2 &= 4\left( d_{i-1}-2d_i+d_{i+1}\right)^2=
\mathcal{O}\left(d_{i-1}^2\right)\,,\\
\beta_3 &= 4 \left( 2d_i-2d_{i+1}\right)^2=\mathcal{O}\left( (d_i-d_{i+1})^2\right)\,,
\end{aligned}
\end{equation*}
resulting in $\beta_2\gg \beta_3$. 
Analogously, if the function $y(x)$ is smooth in $[x_{i-1},x_{i+1}]$
and has {a} discontinuity inside $[x_{i+1},x_{i+2}]$, one obtains
\begin{equation*}
\begin{aligned} 
\beta_2 &=\mathcal{O}\left( (d_i-d_{i-1})^2\right)\,,\\
\beta_3 &=\mathcal{O}\left( d_{i+1}^2\right)\,,
\end{aligned}
\end{equation*}
resulting in $\beta_2\ll \beta_3$.
If the discontinuity point is located in the interval $[x_i,x_{i+1}]$,
both substencils $S_3^{2}$ and $S_3^{3}$ contain the discontinuity.
This seemingly difficult case is actually not problematic, because both polynomials $q_2(x)$ and $q_3(x)$
are essentially monotone inside $[x_i,x_{i+1}]$,
because the over- or undershoots appear in the cells adjacent to the discontinuity \citep{shu1998,shu2009}.
\subsubsection{Smoothness indicators for nonuniform grids}
The first-order four-point derivatives of $y(x)$ at nodes $x_{i-1}$, $x_i$, $x_{i+1}$
and $x_{i+2}$ for nonuniform grids are
approximated by \citep[e.g.,][]{singh2009}\footnote{The formulae in \citet{singh2009} contain some typos, which have been corrected in the present manuscript.}
{\small
\begin{equation}\label{derivatives_nonuniform4}
\begin{aligned}
y'_{i-1} &= -\frac{(2h_{i-1}+h_i)H+h_{i-1}(h_{i-1}+h_i)}{h_{i-1}(h_{i-1}+h_i)H}y_{i-1}
+\frac{(h_{i-1}+h_i)H}{h_{i-1}h_i(h_i+h_{i+1})}y_i\\
&-\frac{h_{i-1}H}{(h_{i-1}+h_i)h_ih_{i+1}}y_{i+1}
+\frac{h_{i-1}(h_{i-1}+h_i)}{(h_i+h_{i+1})h_{i+1}H}y_{i+2}\,,\\ 
y'_i &= -\frac{h_i(h_i+h_{i+1})}{h_{i-1}(h_{i-1}+h_i)H}y_{i-1}
+\frac{h_i(h_i+h_{i+1})-h_{i-1}(2h_i+h_{i+1})}{h_{i-1}h_i(h_i+h_{i+1})}y_i\\
&+\frac{h_{i-1}(h_i+h_{i+1})}{(h_{i-1}+h_i)h_i h_{i+1}}y_{i+1}
-\frac{h_{i-1}h_i}{(h_i+h_{i+1})h_{i+1}H}y_{i+2}\,,\\ 
y'_{i+1} &= \frac{h_ih_{i+1}}{h_{i-1}(h_{i-1}+h_i)H}y_{i-1}
-\frac{h_{i+1}(h_{i-1}+h_i)}{h_{i-1}h_i(h_i+h_{i+1})}y_i\\
&+\frac{(h_{i-1}+2h_i)h_{i+1}-(h_{i-1}+h_i)h_i}{(h_{i-1}+h_i)h_i h_{i+1}}y_{i+1}
+\frac{(h_{i-1}+h_i)h_i}{(h_i+h_{i+1})h_{i+1}H}y_{i+2}\,,\\ 
y'_{i+2} &= -\frac{(h_i+h_{i+1})h_{i+1}}{h_{i-1}(h_{i-1}+h_i)H}y_{i-1}
+\frac{h_{i+1}H}{h_{i-1}h_i(h_i+h_{i+1})}y_i\\
&-\frac{(h_i+h_{i+1})H}{(h_{i-1}+h_i)h_i h_{i+1}}y_{i+1}
+\frac{(2h_{i+1}+h_i)H+h_{i+1}(h_i+h_{i+1})}{(h_i+h_{i+1})h_{i+1}H}y_{i+2}\,, 
\end{aligned}
\end{equation}}\noindent
where $H=h_{i-1}+h_i+h_{i+1}$.
The novel smoothness indicators are constructed as follows
\begin{equation}\label{beta_nonuniform4}
\begin{aligned} 
\beta_2 &= (h_i+h_{i+1})^2\left(\frac{|y'_{i+1}-y'_{i}|}{h_i}-\frac{|y'_{i}-y'_{i-1}|}{h_{i-1}}\right)^2\,,\\
\beta_3 &= (h_{i-1}+h_i)^2\left(\frac{|y'_{i+2}-y'_{i+1}|}{h_{i+1}}-\frac{|y'_{i+1}-y'_{i}|}{h_i}\right)^2\,.
\end{aligned}
\end{equation}
For $h_{i-1} = h_i = h_{i+1}=h$, the finite difference
formulae~\eqref{derivatives_nonuniform4} and~\eqref{beta_nonuniform4}
reduce to uniform grids formulae~\eqref{derivatives_uniform4} and~\eqref{beta_uniform4}, respectively.

In smooth regions {where the second derivative has a constant sign}, the differences of
two consecutive numerical derivatives given by Equation~\eqref{derivatives_nonuniform4} have the same sign, that is,
$$\operatorname{sign} \left( y'_i -y'_{i-1}\right) = \operatorname{sign} \left( y'_{i+1} -y'_i\right) = \operatorname{sign} \left( y'_{i+2} -y'_{i+1}\right)\,.$$
In this case, one shows with some tedious algebra that the indicators~\eqref{beta_nonuniform4}
satisfy $\beta_2=\beta_3$ and, consequently,
\begin{equation*}
\omega_2(x)=\gamma_2(x)\,,\;\text{ and }\;\omega_3(x)=\gamma_3(x)\,.
\end{equation*}
Hence, in smooth regions {where the second derivative has a constant sign}, the nonlinear weights correspond to the optimal weights,
reducing numerical dissipation. 
As above, indicators~\eqref{beta_nonuniform4} allow for over- or undershoots around smooth extrema,
yielding a higher accuracy
near such critical points with respect
{to monotonic interpolations.
As for the uniform case, this feature may have the adverse effect of producing negative interpolated values in strictly
positive function.}

The capability of the smoothness indicators~\eqref{beta_nonuniform4} to detect discontinuities
in the interpolation is already demonstrated for the uniform case in Section~\ref{subsec:uniform4}.
The generalization of the proof to nonuniform grids is particularly cumbersome, but leads to the same conclusion.
However, the indicators~\eqref{beta_nonuniform4} do not guarantee accuracy
near inflection points. 

In conclusion, the smoothness indicators~\eqref{beta_uniform4} and~\eqref{beta_nonuniform4} effectively detect discontinuities and
{ensure} the non-oscillatory property in the interpolation.
Some numerical evidence for the fourth-order WENO interpolations with smoothness indicators~\eqref{beta_uniform4} and~\eqref{beta_nonuniform4}
is displayed in Figures~\ref{fig:exp}-\ref{fig:exp2_not}.
Moreover, an alternative fourth-order WENO interpolation is presented in Appendix~\ref{sec:weno_inter4b}.
\section{Numerical tests}\label{sec:numerical_tests}
For the sake of comparison, 
the third-order ENO,
the third-order WENO \citep[version][]{liu2018}, the cubic Lagrange, the cubic Spline,
{the local monotonic piecewise cubic Hermite \citep{fritsch1984,ibgui2013},}
and the novel fourth-order WENO interpolations
are tested. Four different functions are interpolated on both uniform and nonuniform grids
with different numbers of grid points. The nonuniform grids are randomly generated.
The experimental orders of convergence {on uniform grids} are summarized in Table~\ref{tab:order_accuracy}.
{We note that the monotonic cubic Hermite interpolation is third-order accurate only.
This is due to the derivatives provided by \citet{fritsch1984}, which are second-order accurate on uniform
grids.}
In order to verify the accuracy for smooth cases,
in Figures~\ref{fig:exp} and~\ref{fig:exp_not} we analyze the interpolation of the {exponential function
\begin{equation}\label{exponential}
y(x) = e^{\frac{3}{2}x}\,,\;\text{ with }x\in[-1,1]\,.
\end{equation}
No essential difference between the various interpolations is visible for the homogeneous samplings,
apart from the ENO technique that shows a slightly fragmented interpolation.
Major differences appear for the inhomogeneous samplings, where fourth-order interpolations
perform definitely better and the ENO technique shows a fragmented interpolation.
{Moreover, the monotonic cubic Hermite interpolation proves to be less accurate for
inhomogeneous samplings. This is due to the derivatives provided by \citet{fritsch1984}, which
drop to first-order accuracy on non-uniform grids.}

In order to understand the behaviors around discontinuities,
Figures~\ref{fig:step} and~\ref{fig:step_not} show
the interpolation of the scaled Heaviside function, that is,
\begin{equation}\label{Heaviside}
y(x)= 4H(x)= 
\begin{cases}
0 & \text{ if } x<0 \,, \\ 
4 & \text{ if } x\ge 0\,,
\end{cases}
\;\text{ with }x\in[-1,1]\,.
\end{equation}
{The cubic Lagrange and the cubic Spline interpolations}
reveal Gibbs-like oscillations for homogeneous and inhomogeneous samplings.
{For the homogeneous sampling, such oscillations do not decay in magnitude when the computational grid is refined.}
By contrast, {monotonic cubic Hermite,} ENO and WENO interpolations approximate the discontinuity more sharply
and without oscillations.
In this example, the WENO interpolation algorithm successfully capitalizes on its main idea, {that is},
placing much heavier weights on the candidate substencils in which the discontinuity
is absent.

In order to better understand the different behaviors around smooth critical points (smooth extrema),
Figures~\ref{fig:exp2} and~\ref{fig:exp2_not} analyze the interpolation of the function
\begin{equation}%
\label{exp2}
y(x) = 5(1-e^{-4x^2})\,,\;\text{ with }x\in[-1,1]\,.
\end{equation}
With 16 or more {uniformly} distributed grid points, third-order ENO, the cubic Lagrange, the cubic Spline interpolations, and WENO 4
accurately represent the minimum.
{All these methods require 36 points to accurately represent the extremum for inhomogeneous samplings.}
{By contrast, the accuracy of WENO 3 and of the monotonic cubic Hermite interpolations locally
decreases near the smooth extremum of the interpolated function.}
{We note that the WENO 4 method allows for over- or undershoots around the smooth extrema,
yielding a higher accuracy with respect to monotonic interpolants near such critical points.}

Figures~\ref{fig:sin_step} and~\ref{fig:sin_step_not} show
the interpolation of the function
\begin{equation}\label{sin_step}
y(x)= 
\begin{cases}
2\sin(3x)+4 & \text{ if } x<0 \,, \\ 
2\sin(3x) & \text{ if } x\ge 0\,,
\end{cases}
\;\text{ with }x\in[-1,1]\,,
\end{equation}
which contains both a high-order smooth structure and a discontinuity.
Standard interpolations 
present oscillations for homogeneous and inhomogeneous samplings.
The ENO interpolation is robust and resolves the complex discontinuous function quite well,
{albeit with small overshoots}.
WENO 3 {and the monotonic cubic Hermite interpolation do} not present oscillations,
but {struggle with accurately reproducing the smooth} extrema in the function.
By contrast, WENO 4 shows no oscillations, nor does it produce the spread of the extrema.
%
%
%
\begin{table*}
\caption{Order of accuracy}
\setlength{\tabcolsep}{10pt}\renewcommand{\arraystretch}{1.5}
\centering
\begin{tabular}{ l | c c c c}
\hline
\hline
\emph{Interpolation method}     &\emph{Exponential}     &\emph{Heaviside}       &\emph{Discontinuous sine}   &\emph{Gaussian}\\
\hline 
\hline
ENO                     & 3.017 & 1.072 & 1.073 & 3.033 \\
WENO 3                  & 3.026 & 1.072 & 1.072 & 3.046 \\
Cubic Spline            & 4.043 & 1.024 & 1.025 & 4.046 \\
Cubic                   & 4.043 & 1.029 & 1.029 & 4.044 \\
Monotonic Hermite       & 3.033 & 1.045 & 1.059 & 2.922 \\
WENO 4                  & 4.043 & 1.035 & 1.036 & 4.371 \\\hline
\end{tabular}
\label{tab:order_accuracy}
\end{table*}
%
\section{Application to radiative transfer}\label{sec:RT_applications}
{High accuracy is required in many radiative transfer problems \citep{socas_navarro_2000,trujillo_bueno2003}.
In particular, the upcoming American Daniel K. Inouye Solar Telescope (DKIST) \citep{elmore2014,tritschler2016} and
the planned European Solar Telescope (EST) \citep{matthews2016,matthews2017} have a diffraction limit
near to or surpassing the resolution of the best current 3D radiative-magnetohydrodynamic (R-MHD) simulations of the solar surface \citep{peck2017}, 
making high accuracy critical for comparisons with observations.
Such 3D R-MHD atmospheric models can be highly inhomogeneous and dynamic and present discontinuities.
However, common radiative transfer calculations usually assume smooth variations
in the radiation field and in the atmospheric physical parameters,
and the incidence of discontinuities is often neglected \citep{steiner2016,janett2019}.
In fact, standard numerical schemes and high-order interpolations
tend to misrepresent the nonsmooth regions of
a problem, introducing spurious oscillations near discontinuities \citep[e.g.,][]{richards1991,zhang1997,shu1998}.
Moreover,} many {radiative transfer codes} require grid refinements,
that is, the process of resolving the input data on a finer grid.

The fourth-order WENO interpolation technique presented in Sections~\ref{sec:weno_inter4}
is particularly suitable for scenarios that involve both complex smooth structures and discontinuities
in the physical parameters.
Moreover, it uses the same four-point stencil as the cubic Lagrange and monotonic cubic interpolations.
Therefore, it is fairly straightforward to implement them on already existing codes.
%
%
\subsection{2D and 3D problems}\label{sec:2d_3d_problems}
The only difference between 1D and multidimensional radiative transfer lies
in the mapping of the relevant quantities
along the path of the photons \citep{carlsson2008},
typically through interpolation.
At least for the Cartesian type grids, multidimensional interpolations
can be obtained from 1D procedures. 
In 2D problems, the usual strategy is to discretize the ray by taking its intersections with the segments connecting the grid points,
dealing then with 1D interpolations to recover off-grid points quantities \citep{auer2003}.
In 3D problems, an efficient approach is to discretize the ray by taking its intersections with the individual planes defined by the grid points,
dealing then with 2D interpolations.

{Many radiative transfer codes have the option to use}
linear or bi-linear interpolations from the values of the nearest grid points:
for example, RH by \citet{uitenbroek2001}, MULTI3D by \citet{leenarts2009}, and PORTA by \citet{stepan+trujillo_bueno2013}. 
%
In multiple dimensions, the commonly used short-characteristic
method suffers large numerical diffusion {if} the upwind intensity
is successively {linearly} interpolated in the direction of propagation. {Consequently,}
the solution lacks high-order accuracy \citep{fabianibendicho2003,ibgui2013,peck2017}.
Such {diffusion} errors decrease with increasing order of accuracy of the interpolation
{and a natural choice would be to use} bi-quadratic or bi-cubic interpolations
\citep{kunasz1988,dullemond2013}\footnote{Bi-cubic interpolation can obtained
by using either Lagrange interpolants, cubic convolution algorithm, or cubic splines.}.
However, if the behavior of the quantities to be interpolated is discontinuous or particularly intermittent,
high-order interpolants introduce spurious extrema
that could lead to nonphysical results (e.g., negative values of the source function or of the intensity).
To avoid this, one can opt for monotonic interpolation schemes \citep[e.g.,][]{steffen1990,auer+paletou1994,hayek2010,ibgui2013},
which however degenerate to {a linear interpolation} near smooth extrema, dropping to {second-order accuracy.}
This calls for high-order interpolations that are able to handle discontinuities
and the WENO interpolation presented in Section~\ref{sec:weno_inter4}
can be easily generalized to 2D and 3D Cartesian grids with a direct application
dimension by dimension \citep{zhang2016}.
%
%
\subsection{Multigrid methods}\label{sec:multigrid}
\citet{steiner1991} first implemented a linear multigrid method
for radiative transfer problems.
Later on, \citet{fabianibendicho1997}
applied the nonlinear multigrid method to the multi-level NLTE
radiative transfer problem.
Recently, \citet{stepan+trujillo_bueno2013} generalized the method to the polarized case
and \citet{bjorgen2017} applied a nonlinear multigrid method to realistic
3D R-MHD atmospheric models, which are highly inhomogeneous
and dynamic.

The Jacobi and Gauss-Seidel {iterative schemes} for
3D NLTE radiative transfer problems
quickly smooth out the high-spatial-frequency error,
whereas they slowly decrease the low-frequency error.
To {accelerate the global smoothing of errors},
multigrid methods transfer the problem
to a coarser grid, where such low-frequency errors become
high-frequency errors and iterations on the coarse grid quickly
decrease these errors.
The coarser grid correction is mapped to the finer grid
with an interpolation operator.
{In this step, cubic} interpolations appear to give higher convergence rates~\citep{stepan+trujillo_bueno2013}.
However, when abrupt changes of
atmospheric physical quantities are present (e.g., in 3D R-MHD atmospheric models), such interpolations
introduce spurious oscillations, which
negatively affect the convergence rate of multigrid methods and even induce numerical instabilities,
leading to unphysical negative values of strictly positive physical parameters.
{In such situations, it is safer to use monotonic interpolation, such as the monotonic cubic
Hermite interpolation used in MULTI3D \citep{bjorgen2017}.}
The interpolation operation is {thus} a crucial ingredient of any multigrid method and 
{the fourth-order WENO interpolation might be an
improvement over currently used interpolations.}
%
\subsection{Conversion to optical depth}\label{sec:conv_opt_depth}
{The steady-state version of the radiative transfer equation 
consists in the linear first-order inhomogeneous ODE given by \citep[e.g.,][]{mihalas1978}
\begin{equation}
  \frac{\rm d}{{\rm d} s} I_{\nu}(s) = -\chi_{\nu}(s) I_{\nu}(s) + \epsilon_{\nu}(s)\,,
  \label{eq:scalar_RTE}
\end{equation}
where the spatial coordinate $s$ denotes the position along the ray under consideration, $\nu$ is the frequency,
$I_{\nu}$ is the specific intensity, and
$\chi_{\nu}$ and $\epsilon_{\nu}$ are the absorption and the emission coefficients, respectively. 
In non-local thermodynamic equilibrium (NLTE) conditions, the absorption  and the emission coefficients depend
in a complicated manner on the intensity $I_{\nu}$, so that Equation \eqref{eq:scalar_RTE} is nonlinear
and must be supplemented by additional statistical equilibrium equations
and/or by suitable redistribution matrices for scattering processes.}

The change of coordinates defined by
\begin{equation}
{\rm d}\tau_{\nu}=-\chi_{\nu}(s){\rm d}s\,,
\label{opt_depth}
\end{equation}
transplants Equation~\eqref{eq:scalar_RTE} to the optical depth scale,
namely,
\begin{equation}\label{eq:scalar_RTE_tau}
  \frac{\rm d}{{\rm d} \tau_{\nu}} I_{\nu}(\tau_{\nu}) = I_{\nu}(\tau_{\nu}) - S_{\nu}(\tau_{\nu})\,,
\end{equation}
where $S_{\nu} = \epsilon_{\nu}/\chi_{\nu}$
is the source function.
From Equation~\eqref{opt_depth}, one has
\begin{equation*}
\tau_{\nu} = \tau_{\nu,0} - \int_{s_0}^s \chi_{\nu}(x)\,{\rm d}x\,,
\end{equation*}
and numerical approximations of $\tau_{\nu}(s)$ can be obtained by replacing the integral 
with a numerical quadrature.
\citet{janett2017a} explain that high-order formal solvers require a corresponding high-order quadrature
of the integral above.

A $n$-point weighted quadrature rule is usually stated as 
\begin{equation*}
\int_a^b f(x) {\rm d}x = \sum_{i=1}^n \omega_i f(x_i)\,,
\end{equation*}
and is based on a suitable choice of the nodes $\{x_i\}$ and the weights $\{\omega_i\}$.
In high-order quadratures, one must recover the values of $f(x)$
at off-grid points, {which} typically {requires} high-order interpolations.

As a practical example,
the fourth-order accurate Simpson’s formula
\begin{equation*}
\int_a^b f(x) {\rm d}x = \frac{b-a}{6}\left[f(a)+4f(\tfrac{a+b}{2})+f(b)\right]\,,
\end{equation*}
considers
the node points
$\{a,\tfrac{a+b}{2},b\}\,,$
and their corresponding function values.
The midpoint function value $f(\tfrac{a+b}{2})$
must be recovered through a high-order interpolation, for example, the fourth-order WENO interpolation.
%
%
\subsection{High-order formal solvers}
{When applied to} Equation~\eqref{eq:scalar_RTE} {or Equation~\eqref{eq:scalar_RTE_tau}},
certain high-order
formal solvers,
{such as the RK4 method proposed by \citet{landi_deglinnocenti1976} and
and the third-order pragmatic method by \citet{janett2018b},}
require the absorption $\chi_{\nu}$ and
emission $\epsilon_{\nu}$ coefficients at off-grid points.
In order to maintain high-order accuracy, such quantities must be
obtained through high-order interpolations \citep[e.g.,][for the polarized problem]{janett2017b},
{which} are notoriously oscillatory
in the presence of abrupt changes or discontinuities in the
atmospheric physical quantities.
{In such formal solvers,} spurious oscillations in the interpolation of
$\chi_{\nu}$ and $\epsilon_{\nu}$ can negatively affect their {order of accuracy}.
{The fourth-order WENO interpolation might be befitting in such situations.}
\subsection{Redistribution matrix formalism}\label{sec:red_matrix}
The redistribution matrix formalism is often used for calculating the emissivity in NLTE conditions,
especially when partial frequency redistribution effects are important.
The emissivity is evaluated by
integrating the incident radiation field, multiplied by the redistribution matrix,
over frequency
and direction. 
However, performing the frequency integration is often nontrivial because of the sharp frequency dependence
of the redistribution matrix. 
The quadrature of such integral
requires describing the frequency dependence of the radiation field
between the specific frequency points where its values are known.
Interpolation techniques such as cardinal natural cubic splines
have already been used \citep[see][]{adams1971,belluzzi2014,alsinaballester2017}.
However, such a choice is problematic
if the radiation
field is not sufficiently smooth in the considered frequency interval.
This calls for high-order interpolations capable of dealing with such sharp variations,
such as the fourth-order WENO interpolation.
\section{Conclusions}\label{sec:conclusions}
This paper discusses different ENO and WENO interpolation techniques.
These interpolations are particularly suitable to problems containing both sharp discontinuities and complex smooth structures,
where high-order Lagrange or spline interpolations
are prone to over- and undershoots.

Special attention is paid to the possible applications
of a novel fourth-order WENO interpolation, which is outlined, analyzed, and tested
for both the uniform and nonuniform cases.
This method for computing non-oscillatory fourth-order interpolants
is simple, symmetrical, {and} completely local. 
{Unlike} B\'ezier and monotonic Hermite interpolations,
it does not degenerate
to {second-order accuracy} near smooth extrema of the interpolated function 
and it avoids the use
of conditional ``if$\dots$else'' statements in the algorithm.
Moreover, it uses the same four-point stencil as the cubic Lagrange interpolation.
The implementation on already existing codes is, therefore, fairly straightforward.
Numerical analysis and numerical experiments indicate that the
smoothness indicators given by Equations~\eqref{beta_uniform4} and~\eqref{beta_nonuniform4} guarantee fourth-order accuracy
in smooth regions, while effectively detecting discontinuities and
enabling the non-oscillatory property.
These novel local smoothness indicators yield
less dissipation and better resolution than the canonical ones.
{Such indicators allow for minimal over- or undershoots around smooth extrema,
yielding a higher accuracy near such critical points.
However, we note that this feature may have the adverse effect of producing negative interpolated values in strictly
positive function. Moreover,}
this WENO interpolation is not differentiable at the grid points {(i.e.,} it has kinks).
Even though the fourth-order WENO interpolation may {require}
more CPU time than the cubic Lagrange interpolation 
(depending on the specific algorithm components and computer implementation),
it may still be computationally advantageous for many problems because of its high accuracy.
Moreover, the method can be extended to the interpolation of 2D or 3D data.

{This interpolation technique {might be} particularly suitable in the context of 
the numerical radiative transfer, which must be performed with ever increasing 
spatial and spectral resolution and consequently requires high accuracy in calculations.
However,} the fourth-order WENO interpolation is a general approximation procedure,
{which} can also be used in many other applications beyond radiative transfer,
such as computer vision and image processing. Alternative approximation techniques are available,
such as the so-called HWENO approximations - WENO schemes based on Hermite polynomials \citep{qiu2004,zhu2009}.
\begin{acknowledgements}
The financial support by the Swiss National Science Foundation (SNSF) through grant ID 200021\_159206 is gratefully acknowledged.
{Special thanks are extended to A. Paganini and A. Battaglia for particularly enriching discussions.
The authors are particularly grateful to the anonymous referee for providing valuable comments that helped improving the article.}
\end{acknowledgements}
\begin{figure*}
\includegraphics[width=0.99\textwidth]{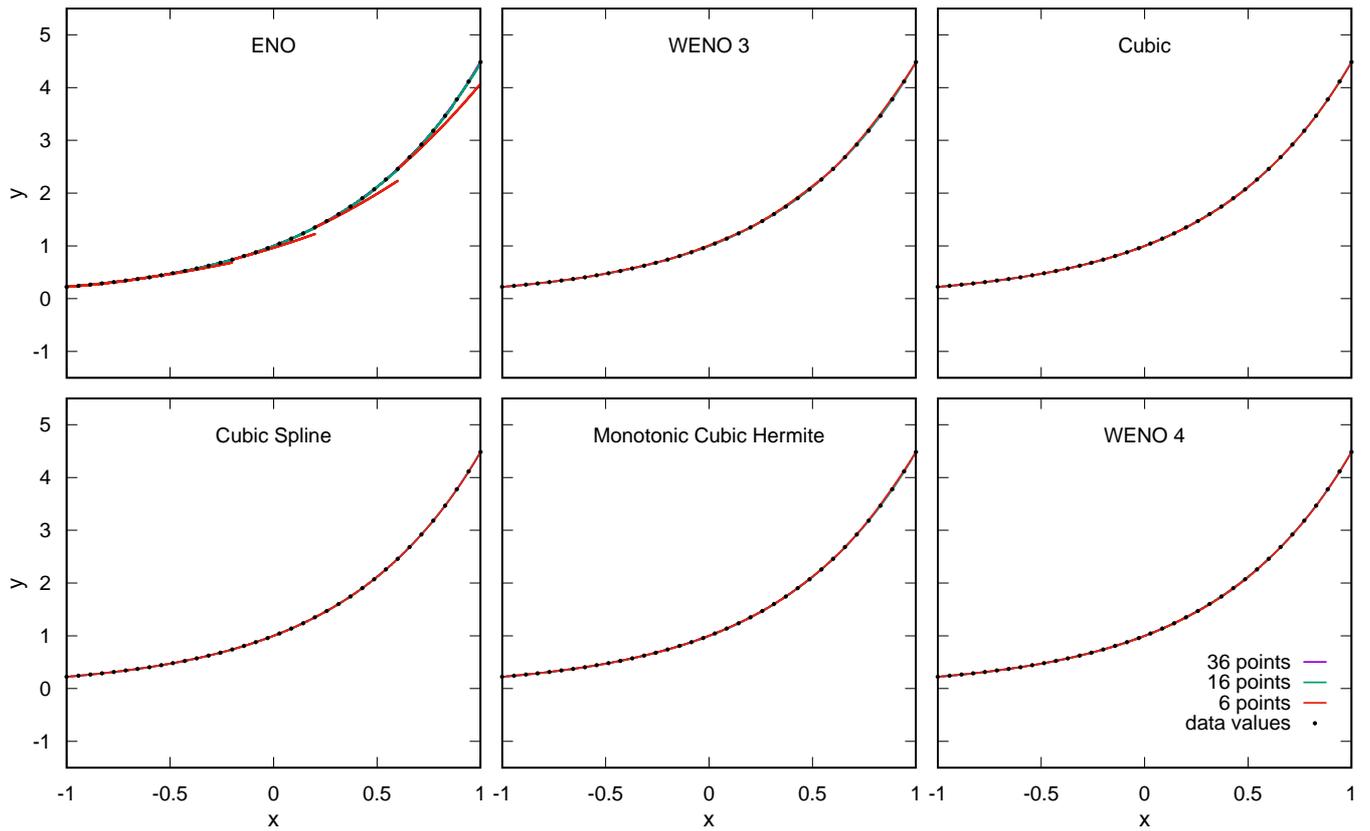}
  \vskip-1.5ex
  \caption{Exponential function~\eqref{exponential},
  is approximated in the interval $x\in[-1,1]$ with
  third-order ENO,
  third-order WENO,
  cubic Lagrange,
  cubic Spline,
  monotonic cubic Hermite,
  and fourth-order WENO interpolations,
  with different homogeneously spaced grid points densities.
  Black dots represent the data values on the 36-point grid.}
  \label{fig:exp}
\end{figure*}
\begin{figure*}
\includegraphics[width=0.99\textwidth]{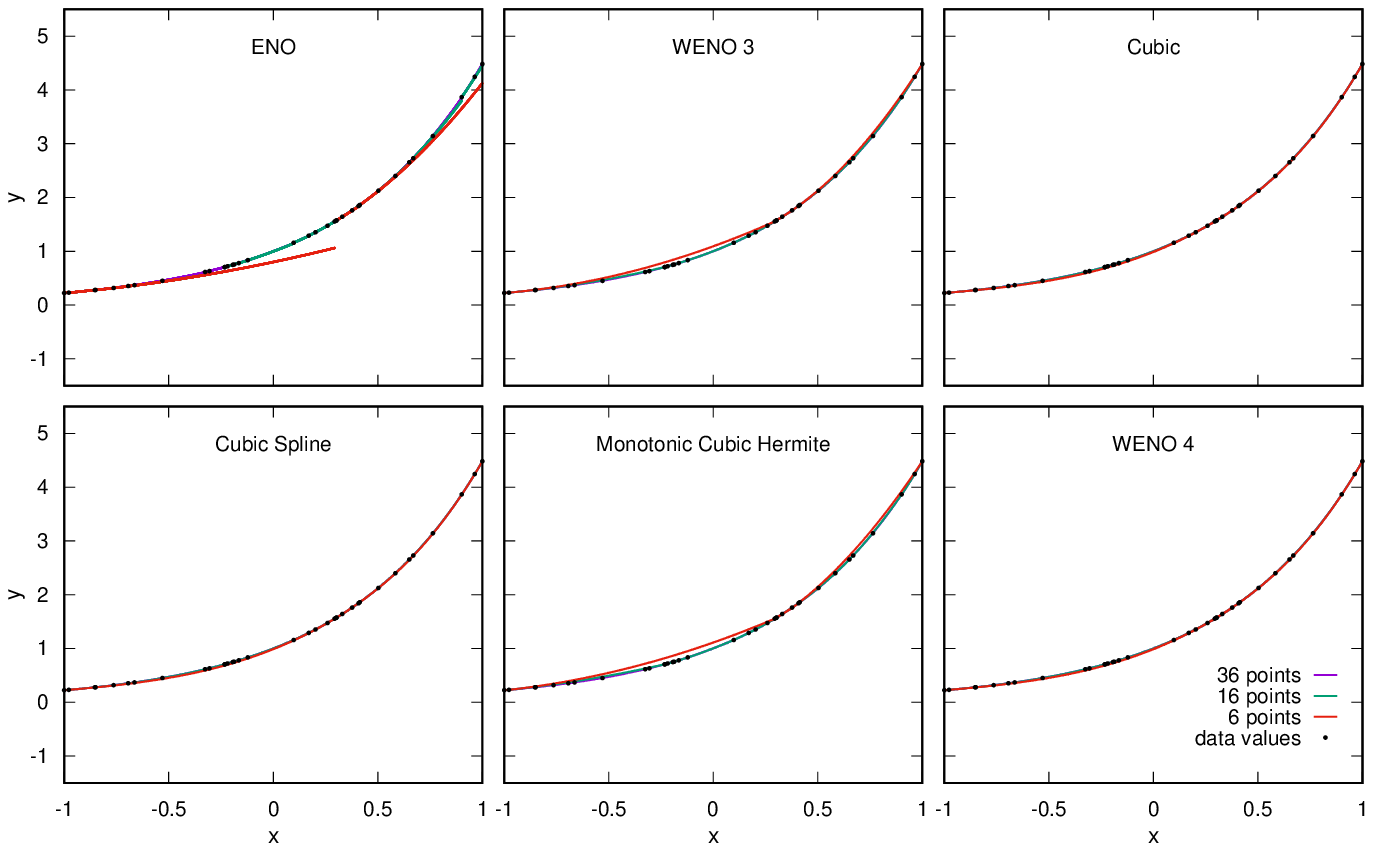}
\vskip-1.5ex
\caption{Same as Figure~\ref{fig:exp}, but for inhomogeneous samplings.}
  \label{fig:exp_not}
\end{figure*}
\begin{figure*}
\includegraphics[width=0.99\textwidth]{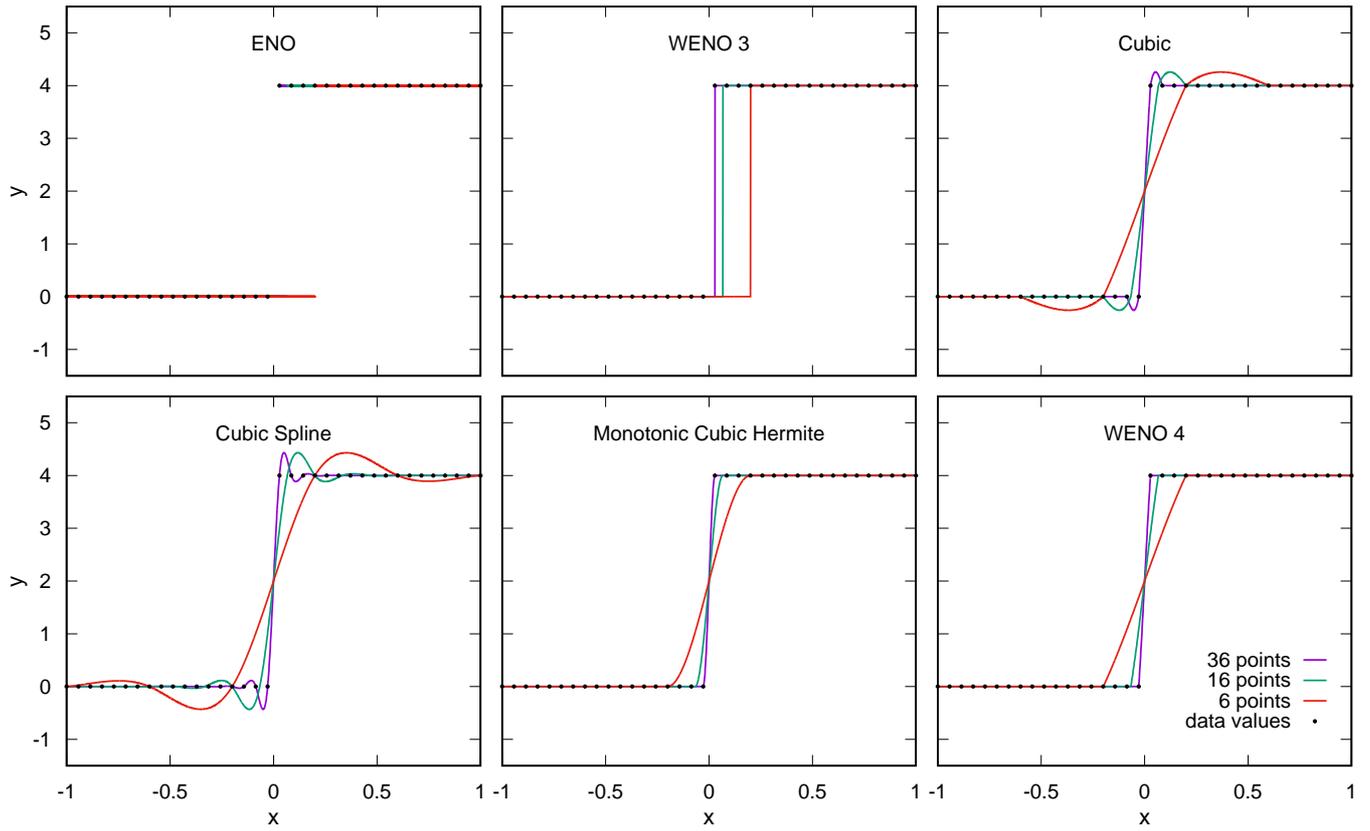}
  \caption{Same as Figure~\ref{fig:exp}, but for the modified Heaviside function given by Equation~\eqref{Heaviside}.}
  \label{fig:step}
\end{figure*}
\begin{figure*}
\includegraphics[width=0.99\textwidth]{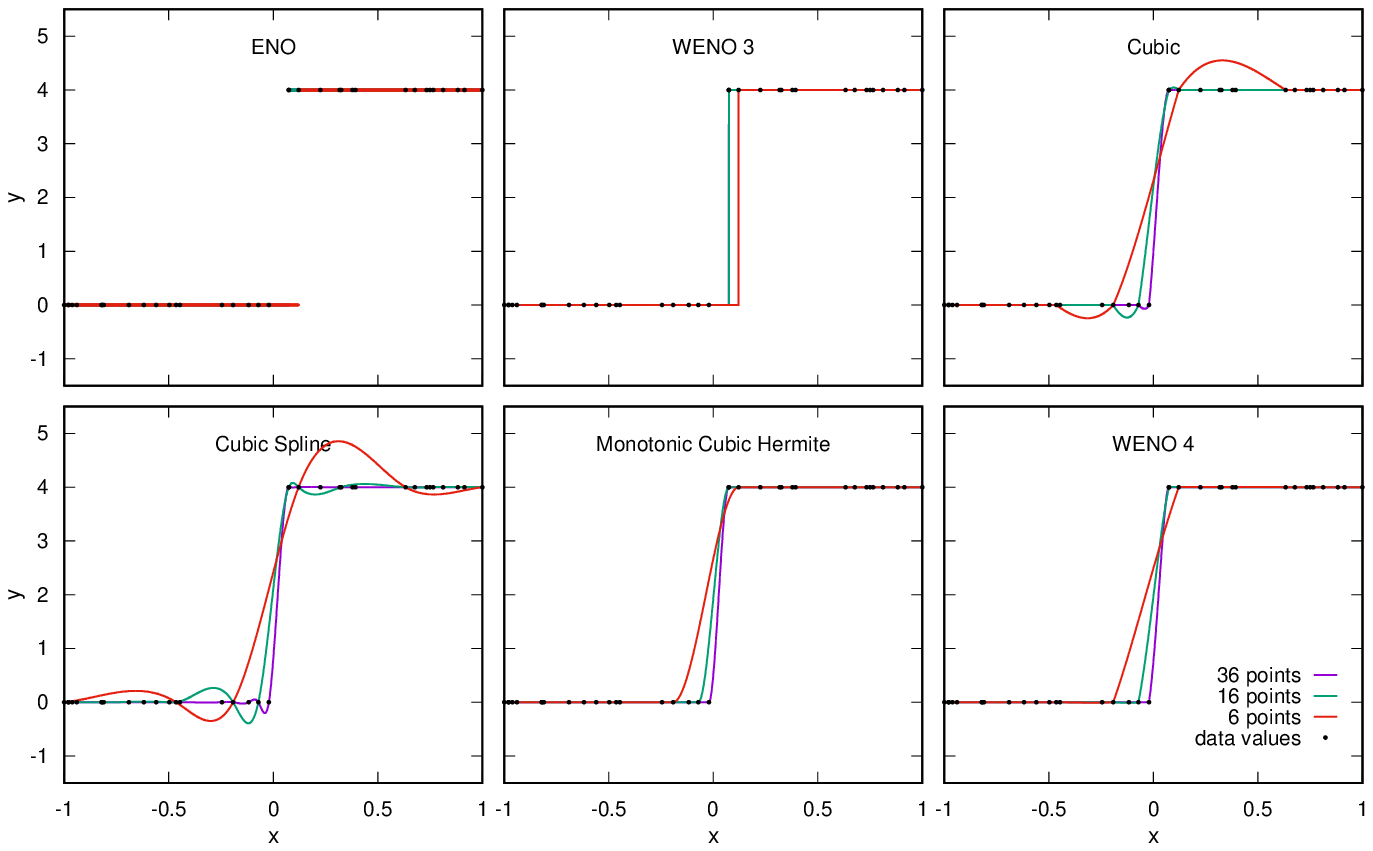}
  \caption{Same as Figure~\ref{fig:step}, but for inhomogeneous samplings.}
  \label{fig:step_not}
\end{figure*}
\begin{figure*}
 \includegraphics[width=0.99\textwidth]{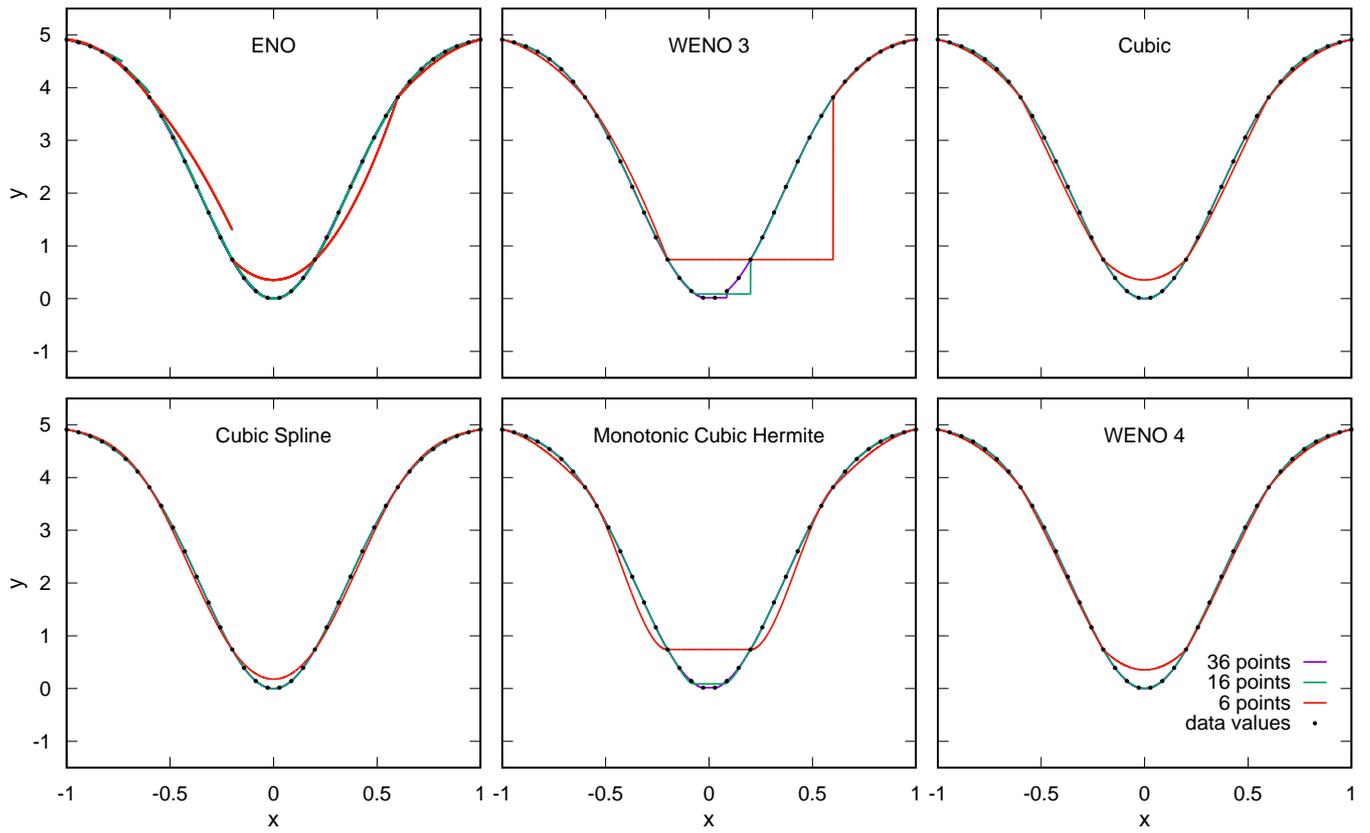}
  \caption{Same as Figure~\ref{fig:exp}, but for the modified Gaussian function given by Equation~\eqref{exp2}.}
  \label{fig:exp2}
\end{figure*}
\begin{figure*}
 \includegraphics[width=0.99\textwidth]{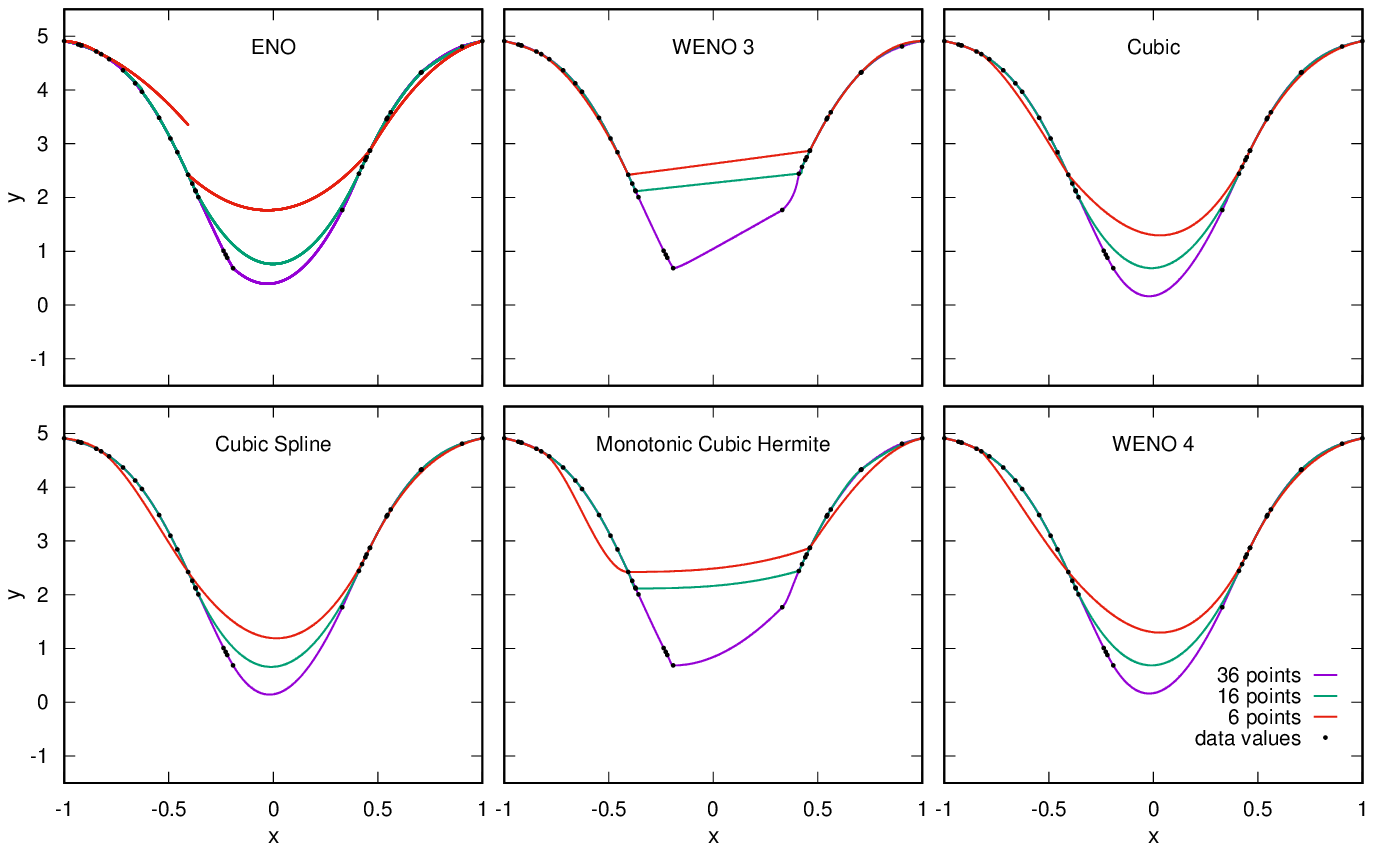}
  \caption{Same as Figure~\ref{fig:exp2}, but for inhomogeneous samplings.}
  \label{fig:exp2_not}
\end{figure*}
\begin{figure*}
 \includegraphics[width=0.99\textwidth]{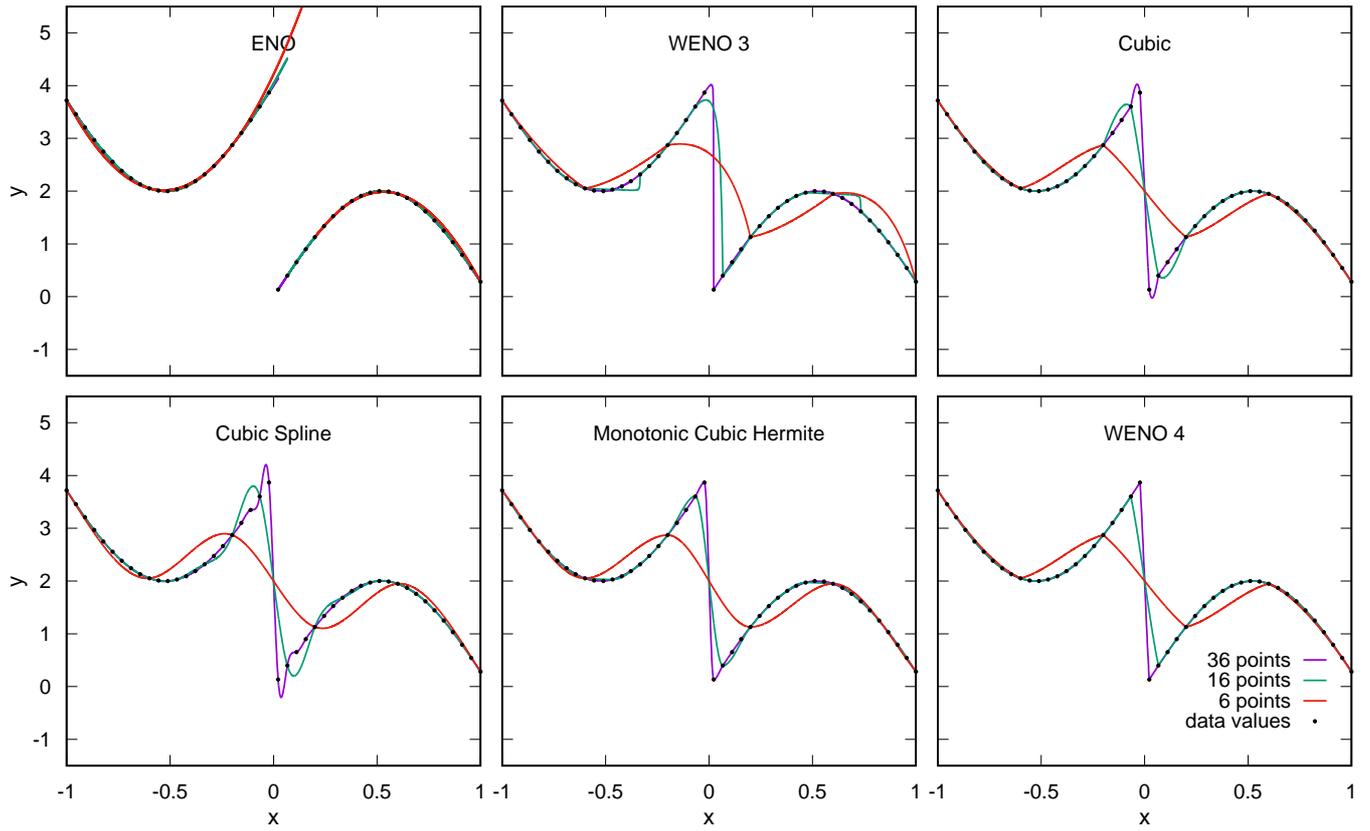}
  \caption{Same as Figure~\ref{fig:exp}, but for the discontinuous sine function given by Equation~\eqref{sin_step}.}
  \label{fig:sin_step}
\end{figure*}
\begin{figure*}
 \includegraphics[width=0.99\textwidth]{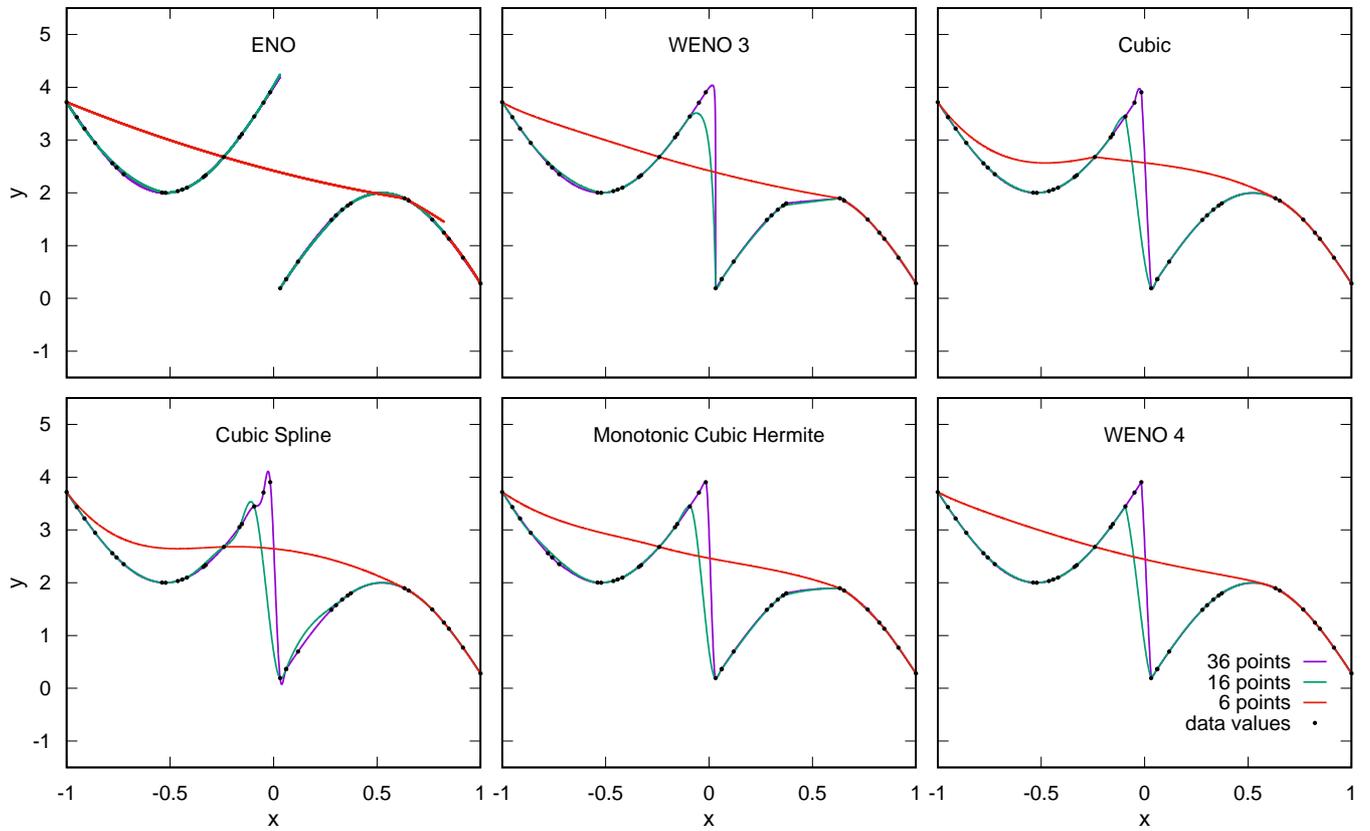}
  \caption{Same as Figure~\ref{fig:sin_step}, but for inhomogeneous samplings.}
  \label{fig:sin_step_not}
\end{figure*}

\appendix

\section{An alternative version of the fourth-order WENO interpolation}\label{sec:weno_inter4b}
In this appendix, an alternative to the fourth-order WENO interpolation presented in Section~\ref{sec:weno_inter4}
is outlined.
The unique cubic Lagrange polynomial $p(x)$ interpolating the four-point large stencil
$$S_4^3 = \{x_{i-1}, x_i , x_{i+1}, x_{i+2}\}\,,$$
can be written as
$$p(x) = \gamma_2(x) q_2(x)+ \gamma_3(x) q_3(x) + \gamma_4(x) q_4(x)\,,$$
in other words, as a linear combination of the linear Lagrange interpolations
\begin{align*}
q_2(x) &= \frac{-y_{i-1}(x-x_i)+y_i(x-x_{i-1})}{x_i-x_{i-1}}\,,\\
q_3(x) &= \frac{-y_{i}(x-x_{i+1})+y_{i+1}(x-x_i)}{x_{i+1}-x_i}\,,\\
q_4(x) &= \frac{-y_{i+1}(x-x_{i+2})+y_{i+2}(x-x_{i+1})}{x_{i+2}-x_{i+1}}\,,
\end{align*}
based, respectively, on the three two-point substencils
$$S_3^{2} = \{x_{i-1}, x_i\}\,,\;S_3^{3} = \{x_i, x_{i+1}\}\,,\; S_3^{4} = \{ x_{i+1}, x_{i+2}\}\,.$$
The linear weights are given by
\begin{align*}
\gamma_2(x) &= \frac{(x-x_{i+1})(x-x_{i+2})}{(x_{i+1}-x_{i-1})(x_{i+2}-x_{i-1})}\,,\\
\gamma_3(x) &= -\frac{(x-x_{i-1})(x-x_{i+2})}{x_{i+2}-x_{i-1}}\left[\frac{1}{x_{i+1}-x_{i-1}}+\frac{1}{x_{i+2}-x_i}\right]\,,\\
\gamma_4(x) &= \frac{(x-x_{i-1})(x-x_i)}{(x_{i+2}-x_{i-1})(x_{i+2}-x_i)}\,,
\end{align*}
and satisfy $\gamma_2(x) + \gamma_3(x) + \gamma_4(x) = 1$.
The linear weights depend just on the grid
geometry and not on the function values.
The fourth-order WENO interpolation in the interval $[x_i,x_{i+1}]$ reads
$$p_i(x) = \omega_2(x) q_2(x)+ \omega_3(x) q_3(x) + \omega_4(x) q_4(x)\,,$$
where the nonlinear weights are defined as
\begin{align*}
\omega_2(x) &= \frac{\alpha_2(x)}{\alpha_2(x)+\alpha_3(x)+\alpha_4(x)}\,,\\
\omega_3(x) &= \frac{\alpha_3(x)}{\alpha_2(x)+\alpha_3(x)+\alpha_4(x)}\,,\\
\omega_4(x) &= \frac{\alpha_4(x)}{\alpha_2(x)+\alpha_3(x)+\alpha_4(x)}\,,
\end{align*}
with
\begin{equation*}
\alpha_2(x) = \frac{\gamma_2(x)}{\epsilon+\beta_2}\,,\;
\alpha_3(x) = \frac{\gamma_3(x)}{\epsilon+\beta_3}\,,\;
\alpha_4(x) = \frac{\gamma_4(x)}{\epsilon+\beta_4}\,,
\end{equation*}
and $\epsilon=10^{-6}$.
\subsection{Smoothness indicators for uniform grids}
In case of uniform grids, a possible choice for the smoothness indicators is
\begin{equation*}
\begin{aligned} 
\beta_2 &= \left(y'_{i}-y'_{i-1}\right)^2\,,\\
\beta_3 &= \left(y'_{i+1}-y'_{i}\right)^2\,,\\
\beta_4 &= \left(y'_{i+2}-y'_{i+1}\right)^2\,,
\end{aligned}
\end{equation*}
where the differences of
the numerical derivatives are given by Equation~\eqref{derivatives_diff_uniform4}.
Unfortunately, this version makes the fourth-order WENO interpolation too dissipative
and generates apparent oscillations.
The search of smoothness indicators that yield
less dissipation and better resolution is not over.
\subsection{Smoothness indicators for nonuniform grids}
In case of nonuniform grids, a possible choice for the smoothness indicators is
\begin{equation*}
\begin{aligned} 
\beta_2 &= \left(\frac{y'_{i}-y'_{i-1}}{h_{i-1}}\right)^2\,,\\
\beta_3 &= \left(\frac{y'_{i+1}-y'_{i}}{h_i}\right)^2\,,\\
\beta_4 &= \left(\frac{y'_{i+2}-y'_{i+1}}{h_{i+1}}\right)^2\,,
\end{aligned}
\end{equation*}
where the numerical derivatives
are given by Equation~\eqref{derivatives_nonuniform4}.
As for the uniform case, this version makes the fourth-order WENO interpolation
too dissipative and
generates apparent oscillations.
The search of smoothness indicators that yield
less dissipation and better resolution is not over.

\bibliographystyle{aa}
\bibliography{bibfile_disc2}

\end{document}